\documentclass[12pt, a4paper]{article}

\usepackage{amsmath}
\usepackage{amsthm}
\usepackage{amssymb}
\usepackage{bbm}
\usepackage{xcolor}
\usepackage{url}

\newtheorem{theorem}{Theorem}

\newtheorem{lemma}[theorem]{Lemma}

\title{Another Variant of the Neoclassical Inequality}

\author{
Stefan Gerhold \\
TU Wien \\
\tt{sgerhold@fam.tuwien.ac.at}
}

\date{\today}

\begin{document}

\maketitle

\begin{abstract}
   The neoclassical inequality is an inequality related to the binomial theorem, involving fractional powers. It naturally gives rise to four other similar inequalities, by
   varying the inequality sign and the range of the fractional parameter.
   Three of these statements are known, and the remaining one was conjectured in the literature,
   based on numerical evidence.  It is proved in the present note. The method we use
   yields short proofs for two of the known related inequalities.
\end{abstract}


The neoclassical inequality~\cite{HaHi10,Ly98} asserts that, for $k\in\mathbb N$
and $0<\alpha<1$,
\begin{equation}\label{eq:nc}
     \alpha \sum_{j=0}^k \binom{\alpha k}{\alpha j} x^{\alpha j} y^{\alpha(k-j)} \leq (x+y)^{\alpha k},
     \quad x,y\geq0,
\end{equation}
with the standard definition
\[
  \binom{\alpha k}{\alpha j} = \frac{\Gamma(\alpha k+1)}{\Gamma(\alpha j+1)\Gamma(\alpha (k-j)+1)}.
\]
See~\cite{GeSi23,HaHi10} for some background and further references. Depending on the
range of $\alpha>0$, maximizing or minimizing the quantity
\[
  F(\alpha,k,x,y) := (x+y)^{-\alpha k}\sum_{j=0}^k \binom{\alpha k}{\alpha j} x^{\alpha j} y^{\alpha(k-j)}, \quad
  k\in\mathbb{N},\ x,y>0,
\]
gives four further inequalities related to~\eqref{eq:nc}, three of which are known. 
The present note completes this quintet, by proving a lower bound for~$F$
in the case $\alpha\geq2$.
As the factor~$\alpha$ on the left hand side of~\eqref{eq:nc} is sharp \cite[p.~468]{HaHi10}, it yields the first line of
\begin{equation}\label{eq:sup}
  \sup\big\{ F(\alpha,k,x,y) : k\in\mathbb N, \ x,y>0 \big\} =
  \begin{cases}
     1/\alpha, & 0<\alpha\leq 1, \\
     1, & \alpha\geq1.
  \end{cases}
\end{equation}
We recall the main result of~\cite{GeSi23}.
\begin{theorem}\label{thm:cnc}
  For $k\in\mathbb N$ and $0<\alpha<1$, 
  \begin{equation}\label{eq:cnc1 xy}
     F(\alpha,k,x,y) \geq 1,
     \quad x,y >0,
  \end{equation}
  and for $\alpha>1$ we have
  \begin{equation}\label{eq:cnc2 xy}
     F(\alpha,k,x,y) \leq 1,
     \quad x,y >0.
  \end{equation}
\end{theorem}
The second line of~\eqref{eq:sup} follows from~\eqref{eq:cnc2 xy}, which is obviously sharp,
as
\[
  \lim_{x\downarrow 0}F(\alpha,k,x,1)=1, \quad k\in\mathbb{N},\ \alpha>0.
\]
The first line of the following statement follows from~\eqref{eq:cnc1 xy}, and the second line from the remark after Theorem~1.3 in~\cite{HaHi10}.
\begin{equation}\label{eq:I}
  \inf\big\{ F(\alpha,k,x,y) : k\in\mathbb N, \ x,y>0 \big\} =
  \begin{cases}
     1, & 0<\alpha\leq 1, \\
     1/\alpha, & 1 \leq\alpha\leq2, \\
     2^{1-\alpha} & \alpha \geq2.
  \end{cases}
\end{equation}
To complete the picture, the following theorem  establishes the last line of~\eqref{eq:I}.
Analogously to~\cite[(1.4)]{GeSi23}, it implies the inequality
\begin{equation}
  E_\alpha((x+y)^\alpha) \leq 2^{\alpha-1} E_\alpha(x^\alpha) E_\alpha(y^\alpha),
  \quad \alpha\geq2,\ x,y \geq0,
\end{equation}
for the Mittag-Leffler function $E_\alpha(x)=\sum_{k=0}^\infty x^k/\Gamma(\alpha k+1)$.

\begin{theorem}\label{thm:main}
    For $k\in\mathbb N$ and $\alpha\geq2$, we have
    \begin{equation}\label{eq:main}
       2^{\alpha-1} \sum_{j=0}^k \binom{\alpha k}{\alpha j} x^{\alpha j}    
              y^{\alpha(k-j)} \geq (x+y)^{\alpha k},
     \quad x,y\geq0.
     \end{equation}
  \end{theorem}
  To see that the factor $2^{\alpha-1}$ is sharp, consider $k=1$, $y=1$
  and $x\to1$.
  In~\cite{GeSi23}, we conjectured Theorem~\ref{thm:main} and noted that it seems
  difficult to prove it with the method from that paper, which employs
  a generalization of the binomial theorem~\cite{HaHi10}. In this note, the theorem is  proved
  with a completely different approach, based on logarithmic convexity w.r.t.~$\alpha$.
\begin{lemma}\label{le:binom}
   For $k>0$ and $0\leq j \leq k$ (not necessarily integers), the function $\alpha \mapsto \binom{\alpha k}{\alpha j}$
   is logarithmically convex on $(0,\infty)$.
\end{lemma}
\begin{proof}
  Define the function~$g$ by $g(0):=0$ and $g(x):=e^{-1/x}/(1-e^{-1/x})$, $x>0$.
  By Lemma~1 in~\cite{KaPr16}, $(\partial^2 /\partial \alpha^2) \log  \binom{\alpha k}{\alpha j}$
  is completely monotone w.r.t.~$\alpha$, and thus non-negative, if 
  \[
    g\Big(\frac{j}{u}\Big) +  g\Big(\frac{k-j}{u}\Big) \leq  g\Big(\frac{k}{u}\Big),
    \quad u>0.
  \]
  As $g$ is convex with $g(0)=0$, it is superadditive~\cite{Br64}, and so
  this criterion is satisfied.
\end{proof}

\begin{proof}[Proof of Theorem~\ref{thm:main}]
 Fix $k\in\mathbb N$, $x,y>0$ and put $p:=x/(x+y)$, $q:=y/(x+y)$, and
 \[
   G(\alpha):=\sum_{j=0}^k \binom{\alpha k}{\alpha j} p^{\alpha j}    
              q^{\alpha(k-j)}, \quad \alpha>0.
 \]
 Clearly, we have $G(1)=(p+q)^k=1$, and
 \begin{align}
   G(2) &= \sum_{j=0}^k \binom{2k}{2j} p^{2j}q^{2(k-j)} 
   = \frac12\big( (p+q)^{2k} + (q-p)^{2k} \big) \notag \\
   &= \frac12 \big( 1 + (q-p)^{2k} \big) \geq \frac12. \label{eq:G est}
 \end{align}
 The statement of Theorem~\ref{thm:main} asserts that $H(\alpha):= \log(2^{\alpha-1}G(\alpha))$ is non-negative for
 $\alpha\geq 2$. As logarithmically convex functions are closed under multiplication (trivial)
 and addition (see~\cite{Ki61}), we deduce from Lemma~\ref{le:binom} that~$G$
 is logarithmically convex on $(0,\infty)$, and thus~$H$ is convex.
 For $\alpha\geq2$, this implies
 \begin{equation}\label{eq:conv}
   H(2) \leq \frac{\alpha-2}{\alpha-1} H(1) + \frac{1}{\alpha-1}H(\alpha)
   = \frac{H(\alpha)}{\alpha-1},
 \end{equation}
 where $H(1)=0$ follows from $G(1)=1$. We have $H(2)=\log(2 G(2))\geq0$ by~\eqref{eq:G est},
 and so~\eqref{eq:conv} gives $H(\alpha)\geq (\alpha-1)H(2)\geq0$ for $\alpha\geq2$.
 Finally, it is clear that~\eqref{eq:main} extends from $x,y>0$ to $x,y\geq 0$ by continuity.
 \end{proof}
While I could not prove~\eqref{eq:nc}
with this method, it does yield a simplified proof of Theorem~\ref{thm:cnc}.
 We keep the notation from the preceding proof. First, consider $0<\alpha<1$.
 Since $1=\frac{1}{2-\alpha}\alpha+\frac{1-\alpha}{2-\alpha}2$, logarithmic convexity of~$G$ yields
\[
 1=G(1) \leq G(\alpha)^{\frac{1}{2-\alpha}} G(2)^{\frac{1-\alpha}{2-\alpha}},
\]
and so $1\leq G(\alpha) G(2)^{1-\alpha}$.
Clearly, $G(2)=\tfrac12(1+(q-p)^{2k})\leq1$.
The latter two estimates and $\alpha\in(0,1)$ imply
\[
  G(\alpha) \geq G(2)^{\alpha-1} \geq1,
\]
which is~\eqref{eq:cnc1 xy}.
We proceed to $\alpha>1$. For $m\in\mathbb{N}$, we have $G(m)\leq 1$, as it is
a partial sum of the binomial expansion of $(p+q)^{mk}=1$. Write
\[
  \alpha = (1-\theta)m+\theta (m+1), \quad m:= \lfloor \alpha \rfloor,\ \theta \in[0,1).
\]
By logarithmic convexity,
\[
  G(\alpha) \leq G(m)^{1-\theta}G(m+1)^{\theta}\leq 1,
\]
and this yields~\eqref{eq:cnc2 xy}.
\bigskip

\noindent{\bf AI disclaimer.} Generative AI (ChatGPT 5.6) was used to explore potential proof strategies and check correctness. No AI-generated material was included without full independent human verification.

\bibliographystyle{siam}
\bibliography{literature}

@article {KaPr16,
    AUTHOR = {Karp, D. B. and Prilepkina, E. G.},
     TITLE = {Completely monotonic gamma ratio and infinitely divisible
              {$H$}-function of {F}ox},
   JOURNAL = {Comput. Methods Funct. Theory},
  FJOURNAL = {Computational Methods and Function Theory},
    VOLUME = {16},
      YEAR = {2016},
    NUMBER = {1},
     PAGES = {135--153},
      ISSN = {1617-9447,2195-3724},
   MRCLASS = {26A48 (33B15 33C60)},
  MRNUMBER = {3460546},
MRREVIEWER = {Valmir\ B.\ Krasniqi},
       DOI = {10.1007/s40315-015-0128-9},
       URL = {https://doi.org/10.1007/s40315-015-0128-9},
}

@article {Ki61,
    AUTHOR = {Kingman, J. F. C.},
     TITLE = {A convexity property of positive matrices},
   JOURNAL = {Quart. J. Math. Oxford Ser. (2)},
  FJOURNAL = {The Quarterly Journal of Mathematics. Oxford. Second Series},
    VOLUME = {12},
      YEAR = {1961},
     PAGES = {283--284},
      ISSN = {0033-5606,1464-3847},
   MRCLASS = {15.60},
  MRNUMBER = {138632},
MRREVIEWER = {J.\ S.\ Frame},
       DOI = {10.1093/qmath/12.1.283},
       URL = {https://doi.org/10.1093/qmath/12.1.283},
}

@article {Ly98,
    AUTHOR = {Lyons, Terry J.},
     TITLE = {Differential equations driven by rough signals},
   JOURNAL = {Rev. Mat. Iberoamericana},
  FJOURNAL = {Revista Matem\'{a}tica Iberoamericana},
    VOLUME = {14},
      YEAR = {1998},
    NUMBER = {2},
     PAGES = {215--310},
      ISSN = {0213-2230},
   MRCLASS = {60H10 (34A99 34F05 37H99)},
  MRNUMBER = {1654527},
MRREVIEWER = {Ben Hambly},
       DOI = {10.4171/RMI/240},
       URL = {https://doi.org/10.4171/RMI/240},
}

@article {HaHi10,
    AUTHOR = {Hara, Keisuke and Hino, Masanori},
     TITLE = {Fractional order {T}aylor's series and the neo-classical
              inequality},
   JOURNAL = {Bull. Lond. Math. Soc.},
  FJOURNAL = {Bulletin of the London Mathematical Society},
    VOLUME = {42},
      YEAR = {2010},
    NUMBER = {3},
     PAGES = {467--477},
      ISSN = {0024-6093},
   MRCLASS = {30B10 (26A33 26D15)},
  MRNUMBER = {2651942},
MRREVIEWER = {Stamatis Koumandos},
       DOI = {10.1112/blms/bdq013},
       URL = {https://doi.org/10.1112/blms/bdq013},
}

@article {Br64,
    AUTHOR = {Bruckner, A. M.},
     TITLE = {Some relationships between locally superadditive functions and
              convex functions},
   JOURNAL = {Proc. Amer. Math. Soc.},
  FJOURNAL = {Proceedings of the American Mathematical Society},
    VOLUME = {15},
      YEAR = {1964},
     PAGES = {61--65},
      ISSN = {0002-9939,1088-6826},
   MRCLASS = {26.52},
  MRNUMBER = {156924},
MRREVIEWER = {E.\ F.\ Beckenbach},
       DOI = {10.2307/2034350},
       URL = {https://doi.org/10.2307/2034350},
}

@article {GeSi23,
    AUTHOR = {Gerhold, Stefan and Simon, Thomas},
     TITLE = {A converse to the neo-classical inequality with an application
              to the {M}ittag-{L}effler function},
   JOURNAL = {Monatsh. Math.},
  FJOURNAL = {Monatshefte f\"ur Mathematik},
    VOLUME = {200},
      YEAR = {2023},
    NUMBER = {3},
     PAGES = {627--645},
      ISSN = {0026-9255,1436-5081},
   MRCLASS = {26D15 (33E12 60G52)},
  MRNUMBER = {4555924},
       DOI = {10.1007/s00605-022-01817-8},
       URL = {https://doi.org/10.1007/s00605-022-01817-8},
}

\end{document}